\documentclass[10pt,reqno]{article}
\usepackage{amsmath,amssymb,amsthm}
\usepackage{esint}
\usepackage{bm}
\usepackage{url}
\usepackage{subfigure}
\usepackage{graphicx,color}
\usepackage{algorithm}
\usepackage{algpseudocode}
\usepackage{algorithmicx}
\usepackage{hyperref}
\usepackage{float}
\usepackage{booktabs,multirow}
\usepackage{hhline}
\usepackage{setspace}

\usepackage{caption}
\usepackage{hyphsubst}
\usepackage{caption}

\usepackage{float}
\usepackage{filecontents}
\usepackage{hyperref}

\hypersetup{colorlinks,citecolor=blue,linkcolor=blue, urlcolor=blue}
\usepackage{pdflscape}
\usepackage{breakurl}

\newcommand{\RR}{\mathbb{R}}
\newcommand{\ds}{\displaystyle}

\DeclareMathAlphabet{\itbf}{OML}{cmm}{b}{it}

\newtheorem{thm}{Theorem}

\numberwithin{equation}{section}

\newcommand{\email}[1]{\protect\href{mailto:#1}{#1}}

\newcommand{\pathfigures}{Figures/}
\graphicspath{{\pathfigures}}

\begin{document}

\title{Comments on \textquotedblleft Generalization of the gradient method with fractional order gradient direction\textquotedblright
}

\author{
Abdul Wahab\footnotemark[1]\, \footnotemark[2]
\and
Shujaat Khan\footnotemark[3]
}
\maketitle
\renewcommand{\thefootnote}{\fnsymbol{footnote}}
\footnotetext[1]{Corresponding Author. E-mail address:  \email{abdul.wahab@sns.nust.edu.pk}.}
\footnotetext[2]{Department of Mathematics, School of Natural Sciences, National University of Sciences and Technology (NUST), Sector H-12, 44000, Islamabad, Pakistan (\email{abdul.wahab@sns.nust.edu.pk)}.}
\footnotetext[3]{Bio-Imaging, Signal Processing, and Learning Lab., Department of Bio and Brain Engineering, Korea Advanced Institute of Science and Technology, 291 Daehak-ro, Yuseong-gu, 34141, Daejeon, South Korea (\email{shujaat@kaist.ac.kr}).}
\renewcommand{\thefootnote}{\arabic{footnote}}

\begin{abstract}

In this paper, a detrimental mathematical mistake is pointed out in the proof of \textit{Theorem 1} presented in the paper\textit{ [Generalization of the gradient method with fractional order gradient direction, J. Franklin Inst., 357 (2020) 2514-2532]}. It is highlighted that the way the authors prove the convergence of the fractional extreme points of a real valued function to its integer order extreme points lacks  correct and valid mathematical argument. Rest of the theorems contained in the paper are mostly announced without any proof relaying on that of Theorem 1. 
\end{abstract}


\noindent {\footnotesize {\bf Key words.} Fractional least mean squares; fractional gradient descent; Fractional calculus; Fractional learning algorithm.}

\section{Introduction}\label{S:Intro}

In \cite{Minima}, a generalized fractional gradient descent scheme is presented and  three fractional least means squares algorithms are introduced.  The main aim of \cite{Minima} is to tackle the problem faced by the fractional order gradient methods in converging to real extreme points. Towards this end, the most important part of the study is the mathematical convergence analysis presented in \cite[Theorem 1]{Minima} that concerns Algorithm 1 \cite[Eq. (11)]{Minima}. Then, two derived algorithms (Algorithm 2 \cite[Eq. (19)]{Minima} and Algorithm 3 \cite[Eq.(29)]{Minima}) are presented and corresponding convergence results are furnished in \cite[Theorems 2 and 3]{Minima} (without any proof by relying on the similarity to that of \cite[Theorem 1]{Minima}).  However, there are some trivial mathematical errors  in the proof of \cite[Theorem 1]{Minima} which are detrimental to the correctness of the entire framework. The main objective of this note is to indicate those mathematical errors. 
\\~\\
\noindent\textbf{Remark:} \textit{The symbols, notations and equation numbers used in this comment are consistent with \cite{Minima}.}

\section{Mathematical Errors}\label{S:Errors}
In order to facilitate ensuing discussion, let us recall \cite[Eq. (4)]{Minima} and \cite[Eq. (11)]{Minima}:
\begin{align}
^C_c\mathcal{D}_x^\alpha f(x) =& \sum_{i=n}^{+\infty}\begin{pmatrix}
\alpha-n\\i-n \end{pmatrix}\frac{f^{(i)}(x)}{\Gamma(i+1-\alpha)} (x-c)^{i-\alpha},
\tag{4}\label{eq4}
\\
x_{k+1} =& x_k -\mu\,\left( ^C_{x_{k-K}}\mathcal{D}_x^\alpha f(x)\right)\Big|_{x=x_k}\qquad (\mu>0,\, K\in\mathbb{Z}_+, \,0<\alpha<1),
\tag{11}\label{eq11}
\end{align}
where 
\begin{align*}
\Gamma(\alpha)=\int_0^{+\infty} e^{-t}t^{\alpha-1} dt \quad (\alpha>0) \quad\text{and}\quad 
\begin{pmatrix}
p\\ q \end{pmatrix} = \frac{\Gamma(p+1)}{\Gamma(q+1)\Gamma(p-q+1)} \quad (p\in\RR, \, q\in\mathbb{N}).
\end{align*}

\subsection{Main Remark}

We show that the proof of \cite[Theorem 1]{Minima} has a detrimental flaw. The statement of \cite[Theorem 1]{Minima} is the following.
\begin{thm}
When the algorithm in \eqref{eq11} is convergent, it will converge to the real extreme point of $f(x)$.
\end{thm}

The method of contradiction is invoked for proof, however, the contraction is obtained through incorrect mathematical argument. We establish our claim below. 

It is assumed that $x^*$ is the real extreme point of $f(x)$ and that the sequence $(x_k)$  converges to a point $X\neq x^*$. Thus, for $0<\varepsilon <|x^*- X|$ there exists $N$ such that $|x_k-X|<\varepsilon < |x^*-X|$ for all $k>N$. Then, by combining \eqref{eq4} and \eqref{eq11}, the following inequality \cite[Eq.(12)]{Minima} is obtained: 
\begin{align}
|x_{k+1}-x_k| = & \mu\left| ^C_{x_{k-K}}\mathcal{D}^\alpha_x f(x)\big|_{x=x_k}\right|
\tag{12a}\label{eq12a}
\\
=& \mu\left|\sum_{i=1}^{+\infty} \begin{pmatrix}\alpha-1\\i-1\end{pmatrix}\frac{f^{(i)}(x_k)}{\Gamma(i+1-\alpha)}\left(x_k-x_{k-K}\right)^{i-\alpha}\right|
\tag{12b}\label{eq12b}
\\
=& \mu\left|\sum_{i=0}^{+\infty} \begin{pmatrix}\alpha-1\\i\end{pmatrix}\frac{f^{(i+1)}(x_k)}{\Gamma(i+2-\alpha)}\left(x_k-x_{k-K}\right)^{i+1-\alpha}\right|
\tag{12c}\label{eq12c}
\\
\geq & \mu\sigma \sum_{i=0}^{+\infty}\left|x_k-x_{k-K}\right|^{i}\left|x_k-x_{k-K}\right|^{1-\alpha}
\tag{12d}\label{eq12d}
\\
= & \mu\sigma \frac{\left|x_k-x_{k-K}\right|^{1-\alpha}}{1-\left|x_k-x_{k-K}\right|}
\tag{12e}\label{eq12e}
\\
\geq & 
d\left|x_k-x_{k-K}\right|^{1-\alpha},
\tag{12f}\label{eq12f}
\end{align}
where 
\begin{align*}
 \sigma:= \sup_{k>N, i\in\mathbb{N}}\begin{pmatrix}\alpha-1\\i\end{pmatrix}\frac{f^{(i+1)}(x_k)}{\Gamma(i+2-\alpha)}, \quad\text{and}\quad d:=\frac{\mu\sigma}{1-\varepsilon}.
\end{align*}
The equations \eqref{eq12a}, \eqref{eq12b}, and \eqref{eq12c} are derived from \eqref{eq11}, \eqref{eq4}, and \eqref{eq12b} (by changing the dummy index $i$), respectively. The inequality \eqref{eq12d} is derived from \eqref{eq12c} by distributing the absolute value over individual terms in the infinite series and then by replacing $\begin{pmatrix}\alpha-1\\i\end{pmatrix}\ds\frac{f^{(i+1)}(x_k)}{\Gamma(i+2-\alpha)}$ by $\sigma$. The argument there to arrive at \eqref{eq12d} is incorrect and the inequality is in the other sense. Indeed, it is trivial to note that by the triangular inequality and Eq. \eqref{eq12c},
\begin{align*}
\mu&\left|\sum_{i=0}^{+\infty} \begin{pmatrix}\alpha-1\\i\end{pmatrix}\frac{f^{(i+1)}(x_k)}{\Gamma(i+2-\alpha)}\left(x_k-x_{k-K}\right)^{i+1-\alpha}\right|
\nonumber 
\\
&\qquad\qquad  \leq \mu\sum_{i=0}^{+\infty} \left|\begin{pmatrix}\alpha-1\\i\end{pmatrix}\frac{f^{(i+1)}(x_k)}{\Gamma(i+2-\alpha)}\right|\,\left|\left(x_k-x_{k-K}\right)^{i+1-\alpha}\right|
\nonumber 
\\
&\qquad\qquad  \leq \mu  \sup_{k>N, i\in\mathbb{N}}\left|\begin{pmatrix}\alpha-1\\i\end{pmatrix}\frac{f^{(i+1)}(x_k)}{\Gamma(i+2-\alpha)}\right|\, \sum_{i=0}^{+\infty} \left|x_k-x_{k-K}\right|^{i}\left|x_k-x_{k-K}\right|^{1-\alpha}.
\end{align*}
Note also that an absolute is missing in the definition of $\sigma$, otherwise it may be negative despite being supremum for a general function $f$, for example, $f(x)=-1/(1-x)$ with $x\in(0,1)$.

In addition to that  \eqref{eq12e} is derived from \eqref{eq12d} by using the sum of infinite geometric series whose convergence is ensured by $|x_k-x_{k-K}|<1$. This tacit assumption is not justified as well for all $k,K\in\mathbb{N}$ even if the series is convergent. 

Finally, inequality \eqref{eq12f} is derived from \eqref{eq12e} with an underlying assumption that $|x_k-x_{k-K}|<\varepsilon$ that also may be true for specific values of $k$ and $K$ but not in general for all $k, K\in\mathbb{N}$ even when the series is convergent.  

In a nutshell, the contradiction derived in \cite[Eqs (13)-(14)]{Minima} using inequality \cite[Eq. (12)]{Minima} is mathematically incorrect.

\subsection{Minor Remarks}
We have following minor comments 
\begin{enumerate}
\item The result announced in \cite[Theorem 2]{Minima} concerning \cite[Algorithm 2]{Minima} is not proved. It is mentioned that \textit{``This theorem can be proved in the similar method like Theorem 1''}.
As we have indicated that the proof of \cite[Theorem 1]{Minima} is flawed, it is unclear whether \cite[Theorem 2]{Minima} is valid or not.

\item A similar remark is valid for \cite[Theorem 3]{Minima}, which is the counterpart of \cite[Theorem 1]{Minima} corresponding to \cite[Algorithm 3]{Minima}.

\item In Eq. \eqref{eq29} of \cite{Minima}, given by 
\begin{align}
x_{k+1}=x_k-\mu\sum_{i=1}^{+\infty} \begin{pmatrix}
\alpha(x)-1\\i-1\end{pmatrix} \frac{f^{(i)}(x_k)}{\Gamma(i+1-\alpha(x))}(x_k-c)^{i-\alpha(x)},
\tag{29}\label{eq29}
\end{align}
it is unclear at which $x$ the exponent $\alpha(x)$ will be evaluated. Moreover, based on the design of $\alpha(x)$ (as discussed in \cite[Fig. 2]{Minima})  $\alpha(x)-1\leq 0$ for all values of $x$ and, thus 
\begin{align*}
\begin{pmatrix}
\alpha(x)-1\\i-1\end{pmatrix} = \frac{\Gamma(\alpha(x))}{\Gamma(i)\Gamma(\alpha(x)-i+1)},
\end{align*}
is undefined for all $i\geq 2$. Indeed, for $0\leq \alpha(x) \leq 1$ and $i\geq 2$, 
\begin{align*}
\alpha(x)-(i-1)\leq \alpha(x)-1\leq 0.
\end{align*}
but the Gamma function $\Gamma(x)$ is defined only for $x\in\RR$ such that $x>0$.
\end{enumerate}

\section*{Conflict of Interest} 
The authors declare that they have no conflict of interest.

\bibliographystyle{plain}

\begin{thebibliography}{99} 



\bibitem{Minima} Y. Wei, Y. Kang, W. Yin, and Y. Wang, Generalization of the gradient method with fractional order gradient direction, {\sl J. Franklin. Inst.}, 357(4): (2020), pp. 2514-2532.

\end{thebibliography}

\end{document}